\documentclass[12pt]{article}
\usepackage{stmaryrd}
\usepackage{mathrsfs}
\usepackage{amsmath}
\usepackage{amsmath,amsthm,amssymb}
\usepackage{amsfonts}
\usepackage{latexsym}
\date{}

\begin{document}
\title{Global dimensions of endomorphism algebras of generator-cogenerators
over $m$-replicated algebras$^\star$}
\author{{\small  Hongbo Lv, Shunhua Zhang}\\
{\small  Department of Mathematics,\ Shandong University,\ Jinan
250100, P. R. China }\\
{\small Dedicated to Professor Shaoxue Liu on the occasion of his
eightieth birthday}}

\pagenumbering{arabic}

\maketitle
\begin{center}
 \begin{minipage}{120mm}
   \small\rm
   {\bf  Abstract} \ \ Let $A$ be a hereditary artin algebra and $A^{(m)}$ be the $m$-replicated
   algebra of $A$. We investigate the possibilities for
   the global dimensions of the endomorphism algebras of generator-cogenerators
   over $A^{(m)}$.


\end{minipage}
\end{center}

\footnote{ $^\star$ Supported by the NSF of China (Grant No.
10771112) and NSF of Shandong Province (Grant No. Y2008A05).}

\footnote{ {\it Email addresses}: lvhongbo356@163.com(H.Lv), \
shzhang@sdu.edu.cn(S.Zhang)}

\section {Introduction}

\vskip 0.2in

Let $\Lambda$ be an artin algebra. We denote by mod $\Lambda$ the
category of all finitely generated right $\Lambda$-modules and by
ind $\Lambda$ a full subcategory of mod $\Lambda$ containing exactly
one representative of each isomorphism class of indecomposable
$\Lambda$-modules. Given a class $\mathcal{C}$ of modules, we denote
by $\rm {add}\ \mathcal{C}$  the subcategory of $\rm{mod}\ \Lambda$
whose objects are the direct summands of finite direct sums of
modules in $\mathcal{C}$ and, if $M$ is a module, we abbreviate
${\rm add}\ \{M\}$ as $ {\rm add }\ M$. The global dimension of
$\Lambda$ is denoted by ${\rm gl.dim}\ \Lambda$ and the
Auslander-Reiten translation of $\Lambda$ by $\tau_\Lambda$. For simplicity, we also write
$\tau_\Lambda$ as $\tau$ if $\Lambda$ is clear.
\vskip 0.2in

A $\Lambda$-module $M$ is called a generator-cogenerator in mod
$\Lambda$ if all indecomposable projective modules and
indecomposable injective modules are in ${\rm add }\ M$.
Let $M$ be a generator-cogenerator in mod $\Lambda$. If ${\rm gl.dim End}_{\Lambda}(M)=d$,
then $M$ is also called a generator-cogenerator with global dimension $d$.

\vskip 0.2in

The endomorphism algebras of generator-cogenerators have attracted a lot
of interest (see for example [T]): these are just the artin algebras
of dominant dimension at least two. The smallest value of the global dimensions of
the endomorphism algebras of such modules was studied by M.Auslander and it was defined
to be the representation dimension of $\Lambda$ provided $\Lambda$ is not semisimple([A]).
If $\Lambda$ is semisimple, its representation dimension is defined to be one. In particular,
M.Auslander proved that an artin algebra $\Lambda$ is representation-finite
if and only if its representation dimension is at most two.

\vskip 0.2in

O.Iyama [I] has shown that the representation dimension
of an artin algebra is always finite, and  R.Rouquier [Rou] has
shown that there is no upper bound for the representation dimensions of artin algebras.
These motivate the investigation on the possibilities for the global
dimensions of the endomorphism algebras of generator-cogenerators.
Recently, V.Dlab and C.M.Ringel described the possibilities for the global
dimensions of the endomorphism algebras of generator-cogenerators in
terms of the cardinalities of the Auslander-Reiten orbits of
indecomposable modules for a hereditary artin algebra (see [DR]).

\vskip 0.2in

Let $A$ be a hereditary artin algebra and $A^{(m)}$ be the $m$-replicated algebra
of $A$. In this paper, we investigate the possibilities for the global
dimensions of the endomorphism algebras of generator-cogenerators in
mod $A^{(m)}$,  and generalize the results in [DR] for $A$ to $A^{(m)}$.
Our main results can be stated as the following.

\vskip 0.2in

{\bf Theorem 1.}\ {\it Let $A^{(m)}$ be the $m$-replicated algebra
of a representation-finite hereditary artin algebra $A$ and $d$ be an
integer with $d\geq 2$. There exists an $A^{(m)}$-module $M$ which
is a generator-cogenerator with global dimension $d$ if and only if
there exists a $\tau_{A^{(m)}}$-orbit of cardinality at least $d$.}

\vskip 0.2in

{\bf Remark.}\ Theorem 1 shows that the possible values for the global
dimension of a generator-cogenerator depend on the maximal length
of the $\tau$-orbits for representation-finite $m$-replicated algebra $A^{(m)}$. That is, for any
integer $b$ with $2\leq b\leq d$, where $d$ is the maximal length
of all $\tau$-orbits, there always exists a generator-cogenerator $M_b$ with
global dimension $b$.

\vskip 0.2in

{\bf Theorem 2.}\ {\it Let $A^{(m)}$ be the $m$-replicated algebra
of a representation-infinite hereditary artin algebra $A$ and let $d$ be
either an integer with $d\geq 3$ or else the symbol $\infty$. Then
there exists an  $A^{(m)}$-module $M$ which is a
generator-cogenerator with  global dimension $d$.}

\vskip 0.2in

{\bf Remark.}\  For every integer $d\geq 3$ or $\infty$,  we will construct a
generator-cogenerator with global dimension $d$  for a given representation-infinite
$m$-replicated algebra $A^{(m)}$. In general,
it is not easy to compute the global dimension of ${\rm End}(M)$
whenever $M$ is a generator-cogenerator,
thus those constructions given in section 4
seem to have an independent interest.

\vskip 0.2in

This paper is arranged as follows. In section 2 we  collect some
definitions and basic facts needed for our research.
Section 3 is devoted to the proof of Theorem 1. In section 4 we construct some
generator-cogenerators in ${\rm mod} \ A^{(m)}$
with special global dimensions and prove Theorem 2.

\vskip 0.2in

\section {Preliminaries}

\vskip 0.2in

Let $\Lambda$ be an artin algebra, $\mathcal{C}$ be a full
subcategory of mod $\Lambda$ and $\varphi :C_M\longrightarrow M$
with $M$ in mod $\Lambda$ and $C_{M}$ in $\mathcal{C}$. The morphism
$\varphi$ is a right $\mathcal{C}$-approximation of $M$ if the
induced  morphism ${\rm Hom}(C,C_{M})\longrightarrow {\rm Hom}(C,M)$
is surjective for any $C$ in $\mathcal{C}$. A minimal right
$\mathcal{C}$-approximation of $M$ is a right
$\mathcal{C}$-approximation which is also a right minimal morphism,
i.e., its restriction to any nonzero direct summand is nonzero. The
subcategory $\mathcal{C}$ is called contravariantly finite if any
module $M$ in mod $\Lambda $ admits a (minimal) right
$\mathcal{C}$-approximation. The notions of (minimal) left
$\mathcal{C}$-approximation and of covariantly finite subcategory
are dually defined. It is well known that, for any module $M$,
the subcategory add $M$ is both a contravariantly finite subcategory
and a covariantly finite subcategory.

\vskip 0.2in

Given a generator $M$, for any $\Lambda$-module $X$, there is a
minimal right ${\rm add}\ M$-approximation of $X$ which is
surjective: $M'\stackrel{f}\longrightarrow X$. Let
$\Omega_M(X)$ be the kernel of $f$. Define inductively $\Omega_M^{i}(X)$ by
$\Omega_M^0(X)=X$, and $\Omega_M^{i+1}(X)=\Omega_M(\Omega_M^{i}(X))$. By definition, the
$M$-$dimension$ $M$-${\rm dim}\ X$ is the minimal value $i$ such
that $\Omega_M^{i}(X)$ belongs to ${\rm add}\ M$ and is $\infty$ if no
$\Omega_M^{i}(X)$ belongs to ${\rm add}\ M$. In particular,
$\Omega_{\Lambda}^i(X)$ is the $i^{\rm th}$ syzygy
and $\Omega^{-i}_{\Lambda}(X)$ is the $i^{\rm th}$ cosyzygy  of
$X$ respectively. We usually write $\Omega_{\Lambda}^i(X)$ (resp. $\Omega^{-i}_{\Lambda}(X)$) as
$\Omega^i(X)$ (resp. $\Omega^{-i}(X)$), and denote by ${\rm pd}\ X$ the projective dimension of $X$. Note that
${\rm pd}\ X= \Lambda$-${\rm dim}\ X$.

\vskip 0.2in

The following lemma is due to Auslander, we refer to [EHIS] for detail.

\vskip 0.2in

{\bf Lemma 2.1.}\ {\it Let $\Lambda$ be an artin
algebra. Let $M$ be a generator-cogenerator and $d\geq 2$. The
global dimension of ${\rm End}\ (M)$ is less or equal to $d$ if and
only if $M$-${\rm dim}\ X\leq d-2$ for any indecomposable
$\Lambda$-module $X$.}

\vskip 0.2in

Let $\Lambda$ be  an artin algebra, and $M, N $ be two
indecomposable $\Lambda$-modules. A path from $M$ to $N$ in ind
$\Lambda$ is a sequence of non-zero morphisms
  $$M=M_0\stackrel{f_1} \longrightarrow M_1\stackrel{f_2} \longrightarrow\cdots
  \stackrel{f_t} \longrightarrow M_t =N$$
  with all $M_i$ in  ind $\Lambda$. Following [Rin], we denote the
  existence of such a path by $M\leq N$. We say that $M$ is a
  predecessor of $N$ (or that $N$ is a successor of $M$).

\vskip 0.2in

More generally, if $S_1$ and $S_2$ are two sets of modules, we
 write $S_1\leq S_2$ if every module in $S_2$ has a predecessor in
 $S_1$, every module in $S_1$ has a successor in
 $S_2$, no  module in $S_2$ has a successor in
 $S_1$ and no  module in $S_1$ has a predecessor in
 $S_2$. The notation $S_1<S_2$ stands for $S_1\leq S_2$
 and $S_1\cap S_2 = \emptyset.$

\vskip 0.2in

 From now on, let $A$ be a hereditary artin algebra and $\mathcal{D}^b({\rm mod}\ A)$
 be the bounded derived category of ${\rm mod}\ A$. The repetitive algebra
 $\hat{A}$ of $A$ was first defined in [HW] and it was proved in [H] that
 $\mathcal{D}^b({\rm mod}\ A)$ is equivalent, as a triangulated
 category, to the stable module category $\underline{{\rm mod}}\ \hat{A}$.

 \vskip 0.2in

The right repetitive algebra $A'$ of $A$ was introduced in [ABM]. Recall from [ABM] that $A'$
is defined as follows,
$$A'=\begin{pmatrix}
             A_{0} &  &  & 0 \\
             Q_{1} & A_{1} &  &  \\
              & Q_{2} & A_{2} &  \\
             0 &  & \ddots & \ddots \\
           \end{pmatrix}
$$
where matrices have only finitely many non-zero coefficients,
$A_{i}=A$ and $Q_{i}= DA$ for all non-negative integers $i$ and $D$ is the standard duality between
${\rm mod}\ A$ and ${\rm mod}\ A^{op}$, all the remaining coefficients are
zero, multiplication is induced from the canonical isomorphisms $A \otimes_A DA\cong\!
DA\cong DA\otimes_A A $ and the zero morphism $DA\otimes_A DA \longrightarrow 0$.

\vskip 0.2in

Recall from [ABST2] that the Auslander-Reiten quiver of $A'$ can be described as follows.
It starts with the Auslander-Reiten quiver of $A_0=A$. Then projective-injective modules
start to appear, such projective-injective module has its socle corresponding to a simple
$A_0$-module, and its top corresponding to a simple $A_1$-module. Next occurs a part
denoted by ${\rm ind}\ A_{01}$ where indecomposables contain at same time simple composition
factors from simple $A_0$-modules, and simple composition factors from simple $A_1$-modules.
When all projective-injective modules whose socle corresponding to simple $A_0$-modules have appeared,
we reach the projective $A_1$-modules and thus the Auslander-Reiten quiver of $A_1$.
The situation then repeats itself.

\vskip 0.2in

{\bf Lemma 2.2.}\ {\it Let $X\in {\rm ind}\ A$, and $\alpha \in {\rm
Hom}_{A'}(\Omega^{-i}_{A'}A, \Omega^{-j}_{A'}X)$.  If $i<j$, then
$\alpha$ factors through a projective-injective $A'$-module. }

\vskip 0.1in

{\bf Proof.} \ The statement follows easily from that
$$\begin{array}{l}
  {\rm\underline{Hom} }_{A'}(\Omega^{-i}_{A'}A,
\Omega^{-j}_{A'}X) \\
  \cong {\rm Hom}_{\mathcal{D}^{b}({\rm mod}\ A
)}(A[i], X[j]) \\
  ={\rm Hom}_{\mathcal{D}^{b}({\rm mod}\ A )}(A,
X[j-i]) \\
  ={\rm Ext}_{A}^{j-i}(A, X) \\
  =0 .
\end{array}$$
 \hfill$\Box$

\vskip 0.2in

By [AI], the $m-$replicated algebra $A^{(m)}$ of $A$ is defined as
the quotient of the right repetitive algebra $A'$, that is,
$$A^{(m)}=\begin{pmatrix}
                  A_{0} &  &  &  & 0 & \\
                  Q_{1} & A_{1} & &  &  &  \\
                   & Q_{2} & A_{2} &  &  &  \\
                 &  &  & \ddots & \ddots &  \\
                  &  0&  &  & Q_{m} & A_{m} \\
                \end{pmatrix}.
$$

If $m=1$, then $A^{(1)}$ is the duplicated algebra of
$A$ (see [ABST1]). According to [AI], we know
that $m+1\leq {\rm gl.dim} \ A^{(m)}\leq 2m+1$.

\vskip 0.2in

Let $\Sigma_0$ be the set of all non-isomorphic indecomposable
projective $A$-modules and set
 $\Sigma_k=\Omega_{A'}^{-k}\,\Sigma_0=\{ \Omega_{A'}^{-k} X \  | \ X\in \Sigma_0 \}$ for $k\geq 0$.
 The next lemma is taken from [ABST2] and will be used in our research.

\vskip 0.2in

{\bf Lemma 2.3.}\ {\it Let $A$ be a hereditary artin algebra. Then

\vskip 0.1in

{\rm (1)}\ The standard embeddings ${\rm ind}\ A_{i}\hookrightarrow
{\rm ind }\ A^{(m)}$ (where $0\leq i\leq m$) and ${\rm ind }\
A^{(m)} \hookrightarrow {\rm ind }\ A'$ are full, exact, preserve
indecomposable modules, almost split sequences and irreducible
morphisms.

\vskip 0.1in

{\rm(2)}\ Let $M$ be an indecomposable $A'$-module which is not
projective and $k\geq 1$. Then the following are equivalent:

\vskip 0.1in

{\rm(a)}\ {\rm pd} $M =k$,

\vskip 0.1in

{\rm(b)}\ $\Sigma_{k-1}< M\leq \Sigma_{k}.$

\vskip 0.1in

 {\rm(3)}\ Let $M$ be an indecomposable $A^{(m)}$-module
which is not in ${\rm ind}\ A= {\rm ind}\ A_0$. Then its projective
cover in {\rm mod} $A^{(m)}$ is projective-injective and coincides
with its projective cover in {\rm mod} $A'$.

\vskip 0.1in

 {\rm(4)}\ Let $M$ be an $A^{(m)}$-module having all
projective-injective indecomposable modules as direct summands. For an
$A^{(m)}$-module $X$, if $X$ has a projective cover which is also projective-injective,
then a minimal right {\rm add} $M$-approximation of $X$ is surjective.}

\vskip 0.2in

In the sequel, we always assume that $A$ is a hereditary artin algebra
and that $A^{(m)}$ is the $m$-replicated algebra of $A$.
It is well known that $A$ and $A^{(m)}$ have the same representation type,
hence $A^{(m)}$ is representation-infinite if and only if $A$ is as well.

\section {The proof of Theorem 1}

This section is devoted to the proof of Theorem 1. We assume in addition
that $A$ is also representation-finite.

\vskip 0.2in

The following result is proved in [DR] for hereditary artin algebras, we
observe that this is also true for $m$-replicated algebras, and the
proof is almost the same, we refer to [DR] for detail. For completeness,
we sketch the proof here.

\vskip 0.2in

{\bf Lemma 3.1.}\ {\it Let $d\geq 2$ be an integer and $M$ be a
generator-cogenerator in  ${\rm mod}\ A^{(m)}$. If any
indecomposable $A^{(m)}$-module $X$ which does not belong to ${\rm
add }\ M$ satisfies that $\tau^{d-1}X=0$, then ${\rm gl.dim}\ {\rm
End}_{A^{(m)}}\ (M)$ is at most $d$.}

\vskip 0.1in

{\bf Proof.} \ It is necessary to prove that $\Omega_M^{d-2}(X)$ is
in ${\rm add}\ M$ for every indecomposable $A^{(m)}$-module $X$ which
does not belong to ${\rm add }\ M$. We follow the idea used in
Proposition 1 of [DR]. Let $Y$ be an indecomposable direct summand of
$\Omega_M^t(X)$. We use induction on $t$ to show that $Y$ is a
predecessor of $\tau^t X$. For $t=d-2$, we see that any
indecomposable direct summand of $\Omega_M^{d-2}(X)$ is a predecessor of
$\tau^{d-2} X$. Since $X$ is not projective and $\tau^{d-1} X=0$, we
know that $\tau^{d-2} X$ is a projective $A$-module, therefore
$\Omega_M^{d-2}(X)$ is also a projective $A$-module, and thus
$\Omega_M^{d-2}(X)$ is in ${\rm add}\ M$. This completes the proof.
\hfill$\Box$

\vskip 0.2in

{\bf Theorem 3.2.}\ {\it Let $A^{(m)}$ be the $m$-replicated algebra
of a representation-finite hereditary artin algebra $A$ and $d$ be an
integer with $d\geq 2$. There exists an $A^{(m)}$- module $M$ which
is a generator-cogenerator with global dimension $d$ if and only if there exists a
$\tau_{A^{(m)}}$-orbit of cardinality at least $d$.}

\vskip 0.1in

{\bf Proof.} \  We assume that there is a $\tau$-orbit of
cardinality at least $d$, thus there exists an indecomposable
$A^{(m)}$-module $Z$ which is not injective such that $\tau^{d-2}Z$ is a projective module.
Let
$$
0\longrightarrow \tau Z\longrightarrow \oplus Y_j\oplus P_Z
\longrightarrow Z\longrightarrow 0
$$
be the Auslander-Reiten sequence ending in $Z$, where all $Y_j$ are
indecomposable and non-projective-injective and $P_Z$ is
projective-injective.

Let $\mathcal{S}=\{\tau^i Z |0\leq i\leq
d-3\}$ and $M$ be the direct sum of all the indecomposables in ${\rm
ind}\ A^{(m)} \backslash \mathcal{S}$. Note that $M$ is a finitely generated
generator-cogenerator since $A^{(m)}$ is representation-finite.

We claim that ${\rm gl.dim}\ {\rm End}_{A^{(m)}} (M)=d.$

In fact, we first show that ${\rm gl.dim}\ {\rm End}_{A^{(m)}} (M)\leq d.$ By
Lemma 2.1, we only need to prove that, for any indecomposable module $X$,
$M$-dim $X\leq d-2.$ If $X$ belongs to ${\rm add}\ M$, then its
$M$-dimension is zero. Thus we only have to consider indecomposable
modules  which do not belong to ${\rm add}\ M$, that is, the modules
$\tau^iZ$ with $0\leq i\leq d-3.$

Let $$0\longrightarrow \tau^{i+1}Z\longrightarrow \oplus
\tau^iY_j\oplus
P_{\tau^iZ}\overset{f_i}{\longrightarrow}\tau^iZ\longrightarrow 0$$
be the Auslander-Reiten sequence ending in $\tau^iZ$, where $0\leq
i\leq d-3$ and $P_{\tau^iZ}$ is a projective-injective module. Since
$A^{(m)}$ is representation-finite, any non-zero non-isomorphism
from an indecomposable module to $\tau^iZ$ can be written as a sum
of compositions of irreducible morphisms. It follows that $f_i$ is
a right ${\rm add}\ M$-approximation and, of course, also minimal.
This shows that $\Omega^i_M(Z)=\tau^iZ$ for $0\leq i\leq d-2$. In
particular, $\Omega^{d-2}_M(Z)=\tau^{d-2} Z$ is projective, thus
belongs to ${\rm add}\ M$. Therefore, $M$-${\rm dim}\ Z\leq d-2$ and consequently,
$M$-${\rm dim}\ \tau^iZ\leq d-2$ for $0\leq i\leq d-3$.

In order to prove that the global dimension of ${\rm End}_{A^{(m)}} (M)$ is
precisely $d$, we show that the $M$-dimension of $Z$ is equal to
$d-2.$ By the construction of $M$, the modules $\Omega^i_M(Z)=\tau^iZ$
for $0\leq i\leq d-3$ are not in ${\rm add }\ M$ and hence $M$-${\rm
dim }\ Z= d-2$.

\vskip 0.1in

For the necessity of Theorem 3.2, we can use Lemma 3.1 and follow
the method used in [DR]. For completeness, we sketch the proof.
Let $M$ be a generator-cogenerator in mod $A^{(m)}$ with
${\rm gl.dim}\ {\rm End}_{A^{(m)}} (M)=d$ and assume that every $\tau$-orbit
is of cardinality at most $d-1$.  Now let $X$ be an indecomposable
$A^{(m)}$-module which is not in ${\rm add} \ M$. Then $\tau^{-1} X$
is non-zero and indecomposable since $X$ is not injective, and we
deduce that $\tau^{d-2} X=0$, since otherwise the sequence of
modules in $\tau^{-1} X, X, \cdots, \tau^{d-2} X$ provides $d$
pairwise non-isomorphic indecomposable modules in a single
$\tau$-orbit.  According to Lemma 3.1, we have that ${\rm gl.dim}\
{\rm End}_{A^{(m)}} (M)\leq d-1$, a contradiction.   This completes the proof
of Theorem 3.2. \hfill$\Box$

\vskip 0.2in

{\bf Remark.}\ We should mention that the proof of Theorem 3.2 heavily depends on the
assumption that $A^{(m)}$ is representation-finite.

\section {Some generator-cogenerators with special global dimensions and the proof of Theorem 2}

\vskip 0.2in

Let $A$ be a hereditary artin algebra  and $A^{(m)}$ be its $m$-replicated algebra.
In this section, we construct some generator-cogenerators in mod $A^{(m)}$
with special global dimensions and give the proof of Theorem 2.

\vskip 0.2in

We denote by $U_{k}$ the direct sum of all the
indecomposable modules in $\Sigma_{k}\bigcap {\rm ind}\ A^{(m)}$ for $k\geq 0$ and
by $P$ the direct sum of all indecomposable projective-injective
$A^{(m)}$-modules. We assume that ${\rm gl.dim}\ A^{(m)}=t$.

\vskip 0.2in

{\bf Proposition 4.1.}\ {\it Let $E_i=A\oplus DA_{m}\oplus P\oplus
\underset{k=i}{\overset{t-1}{\oplus}}U_k$ for $1\leq i\leq t-1$.
Then ${\rm gl.dim} \ {\rm End}_{A^{(m)}}\ (E_i)= i+2$.}

\vskip 0.1in

{\bf Proof.} \ We first show that ${\rm gl.dim} \ {\rm End}_{A^{(m)}}\ (E_i)
\leq i+2$, that is, for each indecomposable $A^{(m)}$-module $X$,
$E_i$-${\rm dim}\ X\leq i.$ If $X $ belongs to ${\rm add}\ E_i$,
then its $E_i$-dimension is zero. Thus we may assume that $X$ is not
in ${\rm add}\ E_i$.

If $\Sigma_0 < X <\Sigma_1$ or $\Sigma_j < X <\Sigma_{j+1}$ with
$i\leq j\leq t-2$ or $\Sigma_{t-1} < X < DA_m$, by using an argument
similar to that in Theorem 3.3 in [LZ], we have that $E_i$-${\rm
dim}\ X\leq 1$ and thus $E_i$-${\rm dim}\ X\leq i$ since $i\geq1$.
For completeness, we give the proof as the following.

Let $\Sigma_0<X<\Sigma_1$, and let
$0 \longrightarrow P_2\longrightarrow P_1 \stackrel{f}\longrightarrow X\longrightarrow 0$
be a minimal projective resolution of $X$.
Since $X\in {\rm ind}\ A$ and ${\rm Hom}_{A^{(m)}}(DA_{m}\oplus
P\oplus \underset{k=i}{\overset{t-1}{\oplus}}U_k, X)=0$, it follows
that $f$ is a minimal right add $E_i$-approximation of $X$ and
$E_i$-${\rm dim}\ X\leq 1$.

Assume that $\Sigma_j<X<\Sigma_{j+1}$ with $i\leq j\leq t-2$
or $\Sigma_{t-1}<X< DA_m$.  It follows that $X$
is not in ind $A$ and there exists an indecomposable $A$-module $Y$
such that $X\cong\Omega^{-j}_{A^{(m)}}Y$ with $i\leq j\leq t-1$.
According to Lemma 2.3(3), $X$ has a projective cover which is projective-injective,
and by Lemma 2.3(4), a minimal right add $(U_j\oplus P)$-approximation of $X$ is surjective.

Consider the following short exact sequence
$$
(\ast)\ \ \ \ \ \ \ 0 \longrightarrow K \longrightarrow M_1 \stackrel{g}\longrightarrow X\longrightarrow 0,
$$
where $g$ is a minimal right add $(U_j\oplus P)$-approximation of
$X$ and $K$ is the kernel of $g$. Let $M=A\oplus DA_{m}\oplus
P\oplus \underset{k=1}{\overset{t-1}{\oplus}}U_k$. Then $g$ is also
a minimal right add $M$-approximation of $X$ by Lemma 2.2 and the
fact that ${\rm Hom}_{A^{(m)}}(DA_{m}\oplus\underset{j<s\leq
t-1}{\oplus}U_s, X)=0$. Since the short exact sequence $(\ast)$ is
not split, $K$ is not projective-injective. Clearly,
$K\leq\Sigma_j$.  Let $N$ be a non-projective-injective
indecomposable direct summand of $K$ and assume that
$\Sigma_l\leq N<\Sigma_{l+1}$, for some $0\leq l\leq j-1$. Then there
exists, by Lemma 2.3(2), an indecomposable $A$-module $Z$ such that
$N\cong\Omega^{-l}_{A^{(m)}}Z$. It follows from Wakamatsu's Lemma
that ${\rm Ext}_{A^{(m)}}^{1}(\Omega^{-(l+1)}_{A^{(m)}}A, K)=0$.
Then we have that $$\begin{array}{rlr}
   0& ={\rm
Ext}_{A^{(m)}}^{1}(\Omega^{-(l+1)}_{A^{(m)}}A, K) \\
   & \cong{\rm
Ext}_{\hat{A}}^{1}(\Omega^{-(l+1)}_{\hat{A}}A, K) \\
   & \cong D {\underline{\rm
Hom}}_{\hat{A}}(\Omega^{-(l+1)}_{\hat{A}}A, \Omega_{\hat{A}}^{-1}K) \\
   & \cong D {\rm
Hom}_{\mathcal{D}^{b}({\rm mod}A)}(A[l+1], K[1])\\
   & \cong D {\rm
Hom}_{\mathcal{D}^{b}({\rm mod}A)}(A, K[-l]).
\end{array}$$
In particular, ${\rm Hom}_{\mathcal{D}^{b}({\rm mod} A)}(A,
N[-l])\cong{\rm Hom}_{A}(A, Z)=0$, and hence $Z=0$ which implies
that $N=0$ and $K\in {\rm add}(U_j\oplus P)$, hence $K\in {\rm add}\ E_i$
and $E_i$-${\rm dim}\ X\leq 1$.

Now, we assume that $\Sigma_l\leq X < \Sigma_{l+1}$ for $1\leq
l\leq i-1$. Let $P(X)$ be a projective cover of $X$. Note that $X$ is not in ${\rm ind}\ A$,
according to Lemma 2.3(3), $P(X)$ is projective-injective, hence we have a short exact sequence
$$
0\longrightarrow \Omega_{A^{(m)}}(X)\longrightarrow
P(X)\stackrel{f_0}\longrightarrow X\longrightarrow 0.
$$
Since ${\rm Hom}_{A^{(m)}}(DA_{m} \oplus
\underset{k=i}{\overset{t-1}{\oplus}}U_k, X)=0$, it follows that
$f_0$ is a minimal right ${\rm add}\ E_i$-approximation of $X$ and
$\Omega_{E_i}(X)=\Omega_{A^{(m)}}(X)$ with $\Sigma_{l-1}\leq \Omega_{E_i}(X)
<\Sigma_l.$ Repeating the step for $\Omega_{E_i}(X)$, we see that
$\Omega_{E_i}^2(X)=\Omega^2_{A^{(m)}}(X)$ and $\Sigma_{l-2}\leq \Omega_{E_i}^2(X)
<\Sigma_{l-1}.$ We repeat this step $l$ times, then
$\Omega_{E_i}^l(X)=\Omega^l_{A^{(m)}}(X)$ and $\Sigma_{0}\leq \Omega_{E_i}^l(X)
<\Sigma_{1}.$ Thus $\Omega_{E_i}^{l+1}(X)$ is just the kernel of a
projective cover of $\Omega_{E_i}^l(X)$ in ${\rm mod}\ A$ and
consequently, $\Omega_{E_i}^{l+1}(X)$ is a projective module and in ${\rm
add}\ E_i$. Therefore, $E_i$-${\rm dim}\ X\leq l+1\leq i.$ This
completes the proof of ${\rm gl.dim} \ {\rm End}_{A^{(m)}}\ (E_i) \leq i+2$ by
Lemma 2.1.

In order to show that ${\rm gl.dim} \ {\rm End}_{A^{(m)}}\ (E_i) $ is exactly
$i+2$, we only need to prove that $E_i$-${\rm dim}\ X =i$ for some
indecomposable module $X$ with $\Sigma_{i-1}< X< \Sigma_{i}$. In
fact, in the argument above, let $l=i-1$, then $\Omega_{E_i}^i(X)$
belongs to  ${\rm  add}\ E_i$ and the modules $\Omega_{E_i}^s(X)$ for
$1\leq s\leq i-1$ are not in ${\rm  add}\ E_i$. It follows that
$E_i$-${\rm dim}\ X =i$. This finishes the proof. \hfill$\Box$

\vskip 0.2in

{\bf Corollary 4.2.}\ {\it Let $d$ be an integer with $3\leq d\leq
t+1$. Then there exists a generator-cogenerator $M$ in ${\rm mod}$
$A^{(m)}$ with global dimension $d$. In particular, the representation
dimension of $A^{(m)}$ is at most 3.}

\vskip 0.2in

From now on, we will assume that $A$ is a hereditary artin algebra
which is representation-infinite. In this case, we observe that $t=2m+1$.

\vskip 0.2in

{\bf Lemma 4.3.}\ {\it Let $A$ be a hereditary artin algebra
which is representation-infinite. Then ${\rm gl.dim}\ A^{(m)}= 2m+1$.}

\vskip 0.1in

{\bf Proof.} \  Let $I$ be an
indecomposable injective $A^{(m)}$-module which is not projective.
By the Auslander-Reiten quiver of $A^{(m)}$, we know that
$I> \Sigma_{2m}$. According to Lemma 2.3(2), we have that ${\rm pd}\ I= 2m+1$,
thus ${\rm gl.dim}\ A^{(m)}= 2m+1$.     \hfill$\Box$

\vskip 0.2in

{\bf Remark.}\ The assumption that $A$ is representation-infinite cannot be omitted.
For example, let $Q$ be the quiver $1\cdot\leftarrow \cdot 2$
and $A=kQ$ be the path algebra of $Q$ over a field $k$.
Then ${\rm gl.dim}\ A^{(1)}= 2$ and ${\rm gl.dim}\ A^{(m)}= m+1$.

\vskip 0.2in

Let $X$ be an $A^{(m)}$-module. If every indecomposable
direct summand of $X$ does not belong to ${\rm ind}\ A$,
we say that $X$ is not in ${\rm mod}\ A$.

\vskip 0.2in

{\bf Lemma 4.4.}\ {\it Let $N=N'\oplus P$ be an $A^{(m)}$-module
with $N'$ being an $A$-module and $P$ the direct sum of all
indecomposable projective-injective $A^{(m)}$-modules. Let $X$ be an
$A^{(m)}$-module which is not in ${\rm mod}\ A$ and $N_X\oplus
P_X\overset{f}{\longrightarrow}X$ be a minimal right ${\rm add}\
N$-approximation with $N_X$ having no projective-injective direct
summands and $P_X$ being projective-injective. Then $P_X$ is a
projective cover of $X$.}

\vskip 0.1in

{\bf Proof.} \ Let $P(X)\overset{g}{\longrightarrow}X$ be a projective cover of $X$.
According to Lemma 2.3 (3) and (4), $P(X)$ is projective-injective and $f$ is surjective.
Then there is a map
$h:\ P(X)\longrightarrow N_X\oplus P_X$ such that $g=fh$. Write
$f=(f_1,f_2)$ with $f_1:\ N_X\longrightarrow X$ and $f_2:\
P_X\longrightarrow X$ and $h=\begin{pmatrix}
                               h_1 \\
                               h_2 \\
                             \end{pmatrix}
$ with $h_1:P(X)\longrightarrow N_X$ and $h_2:P(X)\longrightarrow
P_X.$ Since $N_X$ is in ${\rm add}\ N' $ and also in ${\rm mod}\ A$,
we have that $h_1=0.$ It follows that $g=(f_1,f_2)\begin{pmatrix}
                               0 \\
                               h_2 \\
                             \end{pmatrix}=f_2h_2 $ and thus $f_2$
                             is surjective since $g  $ is
                             surjective. By the minimality of $f $,
                             $P_X$ is isomorphic to $P(X)$. This
                             completes the proof. \hfill$\Box$

\vskip 0.2in

{\bf Lemma 4.5.}\ {\it Let $M$ be a generator-cogenerator in ${\rm
mod}\ A^{(m)}$ such that its non-injective direct summands are all
in ${\rm mod }\ A$. Let $X$ be an indecomposable non-injective
$A^{(m)}$-module. Then $\Omega_M^{2m}(X)$ is an $A$-module. }

\vskip 0.1in

{\bf Proof.} \ Let $\Sigma_0$ be the set of all non-isomorphic indecomposable
projective $A$-modules, and let $\Sigma_k=\Omega_{A'}^{-k}\,\Sigma_0$ for $k\geq
0$. If $X$ belongs to ind $A$, it is trivial
that $\Omega_M^{2m}(X)$ is an $A$-module.

Let $X$ be an indecomposable $A^{(m)}$-module such
that $\Sigma_i\leq X$ for some $i\geq 1$.
Consider the following exact sequence:
$$0\longrightarrow \Omega_M(X)\longrightarrow N_X\oplus P_X
\overset{f}{\longrightarrow}X\longrightarrow 0,$$ where $f$ is a
minimal right ${\rm add}\ M$-approximation and $N_X$ has no
projective-injective direct summands and $P_X$ is
projective-injective. It follows from Lemma 4.4 that $P_X$ is a
projective cover of $X$ and thus $\Omega_M(X)<\Sigma_i.$

Then we only have to show that, for the indecomposables $X$ with
$\Sigma_{2m}\leq X<DA_m$, $\Omega_M^{2m}(X)$ is an $A$-module. By Lemma
4.4 again, $\Omega_M(X)<\Sigma_{2m}$. Write $\Omega_M(X)$ as $K_1\oplus K_2$,
where $K_1$ is an $A$-module and $K_2$ is not in ${\rm mod}\ A$.
There is an exact sequence:
$$0\longrightarrow \Omega_M^2(X)\longrightarrow N'\oplus P'
\overset{g}{\longrightarrow}K_1\oplus K_2\longrightarrow 0,$$ where
$g$ is a minimal right ${\rm add}\ M$-approximation and $N'$ has no
projective-injective direct summands and $P'$ is
projective-injective. By an argument similar to that in the proof of
Lemma 4.4, $P'$ is a projective cover of $K_2$ and thus
$\Omega_M^2(X)<\Sigma_{2m-1}.$ Continuing this procedure, we have that
$\Sigma_0\leq \Omega_M^{2m}(X)< \Sigma_1$ which implies that
$\Omega_M^{2m}(X)$ is an $A$-module. This completes the
proof.\hfill$\Box$

\vskip 0.2in

{\bf Lemma 4.6.}\ {\it Let $M$ be a generator-cogenerator in ${\rm
mod}\ A^{(m)}$ such that its non-injective direct summands are all
in ${\rm mod }\ A$, and let $d$ be an integer with $d\geq 2m+3$. If
$\tau^{d-(2m+2)}N=0$ for any indecomposable non-injective module $N$
in ${\rm add}\ M$, then ${\rm gl.dim}\ {\rm End}_{A^{(m)}} (M)\leq d$. }

\vskip 0.1in

{\bf Proof.} \ According to Lemma 2.1, we have to show that the
$M$-dimension of any indecomposable $A^{(m)}$-module $X$ is at most
$d-2$. If $X $ belongs to ${\rm add}\ M$, then its $M$-dimension is
zero. Thus we only need to consider the indecomposable modules $X$
which do not belong to ${\rm add}\ M$.

Since $X$ is not injective, by Lemma 4.5, $\Omega_M^{2m}(X)$ is an
$A$-module. Let $Y=\Omega_M^{2m}(X)$. Then there is a non-split exact
sequence:
$$0\longrightarrow \Omega_M(Y)\longrightarrow \oplus M_i
\overset{f}{\longrightarrow}Y\longrightarrow 0,$$ where $f$ is a
minimal right ${\rm add}\ M$-approximation of $Y$ and all $M_i $ are
indecomposable modules in ${\rm add}\ M$. Since $Y$ is an
$A$-module, all $M_i$ are non-injective. By our assumption,
$\tau^{d-(2m+2)}M_i=0$ and thus all $M_i$ are preprojective
$A$-modules. This shows that the sequence above is in ${\rm mod}\
A$.

By using induction on $s \geq 1$, one can show that any indecomposable
direct summand $Y'$ of $\Omega_M^s(Y)$ is a predecessor of a non-zero
module of the form $\tau^{s-1}M_i$ for some $i$.  Indeed, we can do it
by an argument similar to that in Proposition 2 in [DR], for completeness,
we include the proof here. If $X'$ is an indecomposable direct summand of $\Omega_M(Y)$,
then $X'$ is a predecessor of some $M_i$. Now, assume the assertion is true
for some $s$. Write $\Omega_M^s(Y)= \oplus Z_j$ with indecomposable modules $Z_j$.
Let $N$ be an indecomposable direct summand of $\Omega_M^{s+1}(Y)=\oplus \Omega_M(Z_j)$.
There is some $j$ such that $N$ is a direct summand of $\Omega_M(Z_j)$.
Note that $\Omega_M(Z_j)\neq 0$, thus $Z_j$ does not belong to ${\rm add}\ M$,
therefore there is a non-split exact sequence
$$
0\longrightarrow \Omega_M(Z_j)\longrightarrow  M
\overset{g}{\longrightarrow} Z_j\longrightarrow 0,
$$
where $g$ is a
minimal right ${\rm add}\ M$-approximation, which means that $N$ is a predecessor
of $\tau Z_j$. By induction, $Z_j$ is a predecessor
of $\tau^{s-1} M_i$ for some $M_i$, it follows that $\tau Z_j$ is a predecessor
of the non-zero $A$-module $\tau^{s} M_i$, since $\tau^{s-1} M_i$ is not projective.
This completes the induction step.

In particular, for $s=d-(2m+2)$, we have that any indecomposable direct summand of
$\Omega_M^{d-(2m+2)}(Y)$ is a predecessor of some $\tau^{d-(2m+3)}M_i$.
By assumption, all the modules $\tau^{d-(2m+3)}M_i$ are projective
$A$-modules, which  implies that
$\Omega_M^{d-(2m+2)}(Y)=\Omega_M^{d-(2m+2)}(\Omega_M^{2m}(X))=\Omega_M^{d-2}(X)$ is also
a projective $A$-module and consequently, belongs to ${\rm add}\ M$.
This completes the proof.\hfill$\Box$

\vskip 0.2in

Let $d$ be an integer with $d\geq 2m+3$. Now we are going to use the
existence of a $\tau$-orbit of cardinality at least $d$ to construct
an $A^{(m)}$-module $M$ which is a generator-cogenerator with global dimension
$d$.

\vskip 0.2in

Since $ A^{(m)}$  is
representation-infinite,  there must be an indecomposable
non-injective $A^{(m)}$-module $Z$ which is also an $A$-module, such that
$\tau^{d-(2m+2)}_{A^{(m)}}Z=\tau^{d-(2m+2)}_{A}Z$ is a
simple and projective $A$-module.

\vskip 0.2in

{\bf Lemma 4.7.}\ {\it  Let $d$ be an integer with $d\geq 2m+3$. Let
$Z$ be an indecomposable non-injective $A^{(m)}$-module such that
$\tau^{d-(2m+2)}Z$ is a simple and projective $A$-module. Let
$$0\longrightarrow \tau Z\longrightarrow \oplus Y_j\longrightarrow
Z\longrightarrow 0$$ be the Auslander-Reiten sequence ending in $Z$,
with indecomposable modules $Y_j$. Let $P$ be the direct sum of all
indecomposable projective-injective $A^{(m)}$-modules. Let
$M=A\oplus DA_m\oplus
\underset{i=0}{\overset{d-(2m+3)}{\oplus}}(\oplus\tau^iY_j)\oplus P$. Then
${\rm gl.dim}\ {\rm End}_{A^{(m)}} (M)=d.$}

\vskip 0.1in

{\bf Proof.} \  We first use Lemma 4.6 to show that ${\rm gl.dim}\
{\rm End}_{A^{(m)}} (M)\leq d.$ Let $N$ be an indecomposable non-injective
module in ${\rm add}\ M$. If $N$ is projective, then $\tau N=0.$ The
assumption of $d\geq 2m+3$ implies that $\tau^{d-(2m+2)}N=0$.
Otherwise, assume that $N=\tau^iY_j$ for some $i$. It is sufficient to show that
$\tau^{d-(2m+2)}Y_j=0.$ Note that $Z $ is actually a preprojective
$A$-module and so all $\tau^iY_j$ are as well. Applying $\tau^{d-(2m+2)}$ to
an irreducible map $Y_j\longrightarrow Z$, we have that either
$\tau^{d-(2m+2)}Y_j=0$ or $\tau^{d-(2m+2)}Y_j$ is a proper
predecessor of $\tau^{d-(2m+2)}Z$. The latter is impossible since
$\tau^{d-(2m+2)}Z$ is simple and projective. Then Lemma 4.6 asserts
that ${\rm gl.dim}\ {\rm End}_{A^{(m)}}\ (M)\leq d$.

In order to show that ${\rm gl.dim}\ {\rm End}_{A^{(m)}} (M)\geq d$,
we only need to find an indecomposable $A^{(m)}$-module $N$ such that
$M$-${\rm dim}\ N \geq d-2$.

For $0\leq i\leq d-(2m+3)$, the Auslander-Reiten
sequence ending in $\tau^iZ$ is of the form
$$
0\rightarrow \tau^{i+1}Z\rightarrow \bigoplus \tau^iY_j
\overset{g_i}{\longrightarrow} \tau^iZ\longrightarrow 0 .
$$
Since all the modules $\tau^iY_j$ with  $0\leq i\leq d-(2m+3)$
belong to ${\rm add}\ M$, we see that $g_i$ is a minimal right
${\rm add}\ M$-approximation.

Note that $\Omega_M^i(Z)=\tau^iZ$ for $0\leq i\leq d-(2m+3)$, this forces that
$M$-${\rm dim}\ Z \geq d-(2m+2)$, since $\tau^{d-(2m+2)}Z$ belongs to ${\rm add}\ M$
and every $\tau^i Z$ does not belong to ${\rm add}\ M$ with $0\leq i<d-(2m+2)$.

Let $N=\Omega^{-2m}_{A^{(m)}}\ Z$. Note that $Z$ is a non-projective preprojective $A$-module,
thus $N$ is an indecomposable $A^{(m)}$-module such that $\Sigma_{2m} <N< DA_m$,
since ${\rm gl. dim}\ A^{(m)}= 2m+1$. It follows from the construction of $M$
that $\Omega^{j}_M(N)=\Omega^j_{A^{(m)}}(N)$
for $1\leq j\leq 2m$, hence $\Omega^{2m}_M(N)=\Omega^{2m}_{A^{(m)}}(N)=Z$.
Therefore  $M$-${\rm dim}\ N = 2m+M$-${\rm dim}\ Z \geq d-2$, and hence
${\rm gl.dim}\ {\rm End}_{A^{(m)}} (M)\geq d$.
This completes the proof.           \hfill$\Box$

\vskip 0.2in

{\bf Lemma 4.8.}\ {\it  Let $N$ be an indecomposable $A$-module
whose endomorphism algebra is a division ring and such that there is
a non-split sequence $0\rightarrow N \stackrel{u}\longrightarrow N'
\stackrel{v}\longrightarrow N \rightarrow 0$.  Let $P$ be the direct
sum of all indecomposable projective-injective $A^{(m)}$-modules and
$M=A\oplus DA_m\oplus P\oplus N'$.  Then $M$ is a
generator-cogenerator in ${\rm mod}\ A^{(m)}$ with ${\rm gl.dim}\
{\rm End}_{A^{(m)}}\ (M)=\infty$.}

\vskip 0.1in

{\bf Proof.} \ One can follow the method used in the proof of
Proposition 3 in [DR] to show that the $M$-dimension of $N$ is
infinite, which implies that ${\rm gl.dim}\ {\rm End}_{A^{(m)}}\ (M)=\infty$.
For the convenience of readers, we streamline the proof.

Let $p: P(N)\rightarrow N$ be a projective cover of $N$. Then
$[v,p]: N'\oplus P(N)\rightarrow N$ is a right ${\rm add}\ M$-approximation
(may-be not minimal) since any map $f:
N'\rightarrow N$ factors through $v$, and any map $P \rightarrow N$ with $P$
projective factors through $p$, and ${\rm Hom}_{A^{(m)}} (I,N)=0$ for any injective
module $I$.

Note that $p$ itself is not a right ${\rm add}\ M$-approximation,
since the map $v$ cannot be factored through a projective module.
Thus, a minimal right ${\rm add}\ M$-approximation of $N$ is of the
form $[v,p']: N'\oplus P'\rightarrow N$ with $P'$ projective.

The kernel of $[v,p']$ is isomorphic to $N\oplus P''$ with $P''$
projective, thus $\Omega_M(N)= N\oplus P''$. Inductively, we see that
$N$ is a direct summand of $\Omega_M^t(N)$ for all $t\geq 0$, this
shows that the $M$-dimension of $N$ is not finite.
The proof is completed.       \hfill$\Box$

\vskip 0.2in

{\bf Theorem 4.9.}\ {\it Let $A^{(m)}$ be the $m$-replicated
algebra of a representation-infinite hereditary artin algebra $A$ and let
$d$ be either an integer with $d\geq 3$ or else the symbol $\infty$.
Then there exists an $A^{(m)}$-module $M$ which is a
generator-cogenerator with global dimension $d$.}

\vskip 0.1in

{\bf Proof.} \  If $d$ is an integer with $3\leq d\leq 2m+2$, according to Corollary 4.2,
there exists a generator-cogenerator $M$ in mod $A^{(m)}$
with ${\rm gl.dim}\ {\rm End}_{A^{(m)}} (M)=d$.

If $d$ is an integer with $d\geq 2m+3$, then the consequence follows
from Lemma 4.7.

Finally, let $d=\infty$. Since $A^{(m)}$ is representation-infinite,
it follows that $A$ is representation-infinite, and thus there must
exist an indecomposable $A$-module $N$ whose endomorphism
algebra is a division ring and such that ${\rm
Ext}_{A}^{1}(N,N)={\rm Ext}_{A^{(m)}}^{1}(N,N) \neq 0$. Then
according to Lemma 4.8, there exists a generator-cogenerator $M$ in
mod $A^{(m)}$ such that ${\rm gl.dim}\ {\rm End}_{A^{(m)}} (M)=\infty$. This
completes the proof. \hfill$\Box$

\vskip 0.2in

{\bf Acknowledgements.} \ We thank the referee for his or her
valuable suggestions and comments.

\begin{description}

\item{[A]}\ M.Auslander. Representation Dimension of Artin
Algebras. Queen Mary College Mathematics Notes, Queen Mary College,
London, 1971.

\item{[ABM]}\ I.Assem, A.Beligiannis, N.Marmaridis, Right
triangulated categories with right semi-equivalences. {\it CMS
Conference Proceedings} 24(1998), 17-37.

\item{[ABST1]}\ I.Assem, T.Br$\ddot{\rm u}$stle, R.Schiffer,
G.Todorov,  Cluster categories and duplicated algebras. {\it J.
Algebra} 305(2006), 548-561.

\item{[ABST2]}\ I.Assem, T.Br$\ddot{\rm u}$stle, R.Schiffer,
G.Todorov,  $m$-cluster categories and $m$-replicated algebras. {\it
Journal of pure and applied Algebra} 212(2008), 884-901.

\item{[AI]}\ I.Assem, Y.Iwanaga,  On a class of
representation-finite QF-3 algebras. {\it Tsukuba J. Math. }
11(1987), 199-210.

\item{[DR]}\ V.Dlab, C.M.Ringel, The global dimension of the
endomorphism ring of a generator-cogenerator for a hereditary artin
algebra. {\it Mathematical Reports of the Academy of Science of the Royal
Society of Canada}. 30(3)(2008), 89-96.

\item{[EHIS]}\ K.Erdmann, T.Holm, O.Iyama, J.Schroer, Radical embedding and
representation dimension, Adv.Math. 185(2004), 159-177.

\item{[H]}\ D.Happel, Triangulated categories in the representation
theory of finite dimensional algebras. {\it Lecture Notes series
119.} Cambridge Univ. Press, 1988.

\item{[HW]}\ D.Hughes, J.Waschb$\ddot{\rm u}$sch,  Trivial extensions
of tilted algebras. {\it Proc.London. Math.Soc.} 46(1983),  347-364.

\item{[LZ]}\ H.Lv, S.Zhang, Representation dimension of
$m$-replicated algebras. Preprint, arXiv:math.RT/0810.5185, 2008.

\item{[I]}\ O.Iyama, Finiteness of representation dimension.
{\it Proc.Amer.Math.Soc.}, 131(2003), 1011-1014.

\item{[Rin]}\ C.M.Ringel, Tame algebras and integral quadratic forms.
{\it Lecture Notes in Math. 1099.} Springer Verlag, 1984.

\item{[Rou]}\ R. Rouquier,  Dimensions of triangulated categories.
{\it J. K-theory: K-theory and its Applications to Algebra, Geometry, and Topology. } 1(2008), 193-256.

\item{[T]}\ H.Tachikawa, {\it Lectures on Quasi-Frobenius Rings.}
Springer Lecture Notes in Mathematics 351 (1973).

\end{description}

\end{document}